\documentstyle[amstex,amssymb,12pt]{article}
\input cyracc.def
\newfam\cyifam
\font\tencyi=wncyi10
\font\sevencyi=wncyi7
\font\fivecyi=wncyi5

\textfont\cyifam=\tencyi \scriptfont\cyifam=\sevencyi
  \scriptscriptfont\cyifam=\fivecyi

\newfam\cyrfam
\font\tencyr=wncyr10
\font\sevencyr=wncyr7
\font\fivecyr=wncyr5
\def\cyr{\fam\cyrfam\tencyr\cyracc}
\textfont\cyrfam=\tencyr \scriptfont\cyrfam=\sevencyr
  \scriptscriptfont\cyrfam=\fivecyr

\newcommand{\lon}{\longrightarrow}
\newcommand{\ad}{{\mathrm ad}}
\newcommand{\rar}{\rightarrow}

\newcommand{\Proof}{{\bf Proof}.\, }

\newcommand{\Def}{\mathsf{Def}}
\newcommand{\Df}{\mbox{\cyr Def}}

\newcommand{\Z}{{\Bbb Z}}
\newcommand{\p}{{\partial}}

\newcommand{\C}{{\Bbb C}}
\newcommand{\R}{{\Bbb R}}
\newcommand{\ot}{\otimes}
\newcommand{\tl}{\tilde}
\newcommand{\tlk}{\tilde{\kappa}}

\newcommand{\Beq}{\begin{equation}}
\newcommand{\Eeq}{\end{equation}}
\newcommand{\Beqr}{\begin{eqnarray}}
\newcommand{\Eeqr}{\end{eqnarray}}
\newcommand{\Beqrn}{\begin{eqnarray*}}
\newcommand{\Eeqrn}{\end{eqnarray*}}
\newcommand{\Ba}{\begin{array}}
\newcommand{\Ea}{\end{array}}
\newcommand{\Bi}{\begin{itemize}}
\newcommand{\Ei}{\end{itemize}}
\newcommand{\Bc}{\begin{center}}
\newcommand{\Ec}{\end{center}}
\newcommand{\fg}{{\frak g}}

\newcommand{\fm}{{\frak m}}
\newcommand{\f}{{\cal O}}

\newcommand{\cM}{{\cal M}}

\newcommand{\cB}{{\cal B}}

\newcommand{\al}{\alpha}

\newcommand{\ga}{\gamma}
\newcommand{\ka}{\kappa}
\newcommand{\Ga}{\Gamma}

\newcommand{\tlo}{\tilde{0}}
\newcommand{\tla}{\tilde{a}}
\newcommand{\tln}{\tilde{1}}

\newcommand{\bp}{\bar{\partial}}

\newcommand{\Ker}{{\mathsf Ker}\, }
\newcommand{\Img}{{\mathsf Im}\, }

\newcommand{\cT}{{\cal T}}

\newcommand{\sip}{\smallskip}
\newcommand{\bip}{\bigskip}

\newcommand{\vse}{\vspace{8 mm}}

\begin{document}

\title{Semi-infinite $A$-variations of Hodge structure\\ over extended K\"ahler cone}
\author{S.A.\ Merkulov}
\date{}
\maketitle



\sloppy
\begin{center}
{\bf \S 1.  Introduction}
\end{center}

\bip

According to Kontsevich's Homological Mirror Conjecture \cite{K1},
a {\em mirror pair}, $X$ and $\hat{X}$,
  of
Calabi-Yau manifolds  has two associated
$A_{\infty}$-categories, the derived category of coherent sheaves
on $X$ and the Fukaya category of $\hat{X}$, equivalent. In
particular, the moduli spaces of $A_{\infty}$-deformations of these
two categories must be isomorphic implying
$$
H^*(X, \wedge^*T_X) = H^*(\hat{X},\C).
$$
Another expected corollary is the equivalence of two Frobenius manifold
structures, the first one is generated on $H^*(X,\wedge^* T_X)$
by the
periods of semi-infinite variations of Hodge structure on X \cite{Ba,Ba2},
and the second one is generated on $H^*(\hat{X},\C)$  by the Gromov-Witten
invariants.

\sip

The l.h.s.\ in the above equality can be identified with the
tangent space at $X$ to the extended moduli space, $\cM_{\mathrm compl}$, of complex
structures \cite{BK}. This moduli space is the base of
semi-infinite $B$-variations, ${\mathsf VHS}^B(X)$, of the standard Hodge
structure in $H^*(X,\C)$ \cite{Ba}. Moreover, it was shown in
\cite{Ba}
how to construct a family of Frobenius manifold structures,
  $\{\Phi^W_{\mathrm compl}(X)\}$,
  on $\cM_{\mathrm compl}$
  parameterized by isotropic increasing filtraions, $W$, in
  the de Rham cohomology $H^*(X,\C)$ which are complementary to
  the standard decreasing Hodge filtration. Presumably,
  the compactification  $\overline{\cM}_{\mathrm compl}$ contains a point with
maximal unipotent monodromy, and the associated limiting weight
filtration $W_0$ gives rise, via Barannikov's semi-infinite variations of Hodge structure,
 to the solution, $\Phi^{W_0}_{\mathrm compl}(X)$,  of the WDVV equations
 which coincides precisely with the potential, $\Phi_{GW}(\hat{X})$, built out of the
 Gromov-Witten invariants on the mirror side. This has been
 checked for complete Calabi-Yau intersections in \cite{Ba0}.

 \sip

 It is widely believed that $\Phi_{GW}(\hat{X})$ can itself be
 reconstructed from $A$-model variations of  Hodge
 structure (see \cite{CF,CK,Mo} for the small quantum cohomology
 case).
In this paper we propose a symplectic version of the Barannikov's
construction which, presumably, provides a correct framework for
extending the results of \cite{CF,CK,Mo} to the full quantum
cohomology group.
 We study semi-infinite $A$-variations,
${\mathsf VHS}^A(\hat{X})$, of Hodge
structure over the extended moduli space, $\cM_{\mathrm sympl}$,  of K\"ahler
forms on the mirror partner $\hat{X}$, and then use Barannikov's technique \cite{Ba}
to build out of ${\mathsf VHS}^A(\hat{X})$ a family
 of solutions of the WDVV equations,
$\{\Phi_{\mathrm sympl}^{W'}(\hat{X})\}$, parameterized by isotropic increasing
filtraions, $W'$, in $H^*(X,\wedge^* T_X)$ which are complementary to
  the standard Hodge type filtration in $\oplus_{i,j}H^i(X,
  \wedge^j T_X)$. The tangent space to $\cM_{\mathrm sympl}$ is
   $H^*(\hat{X}, C)$ which is  the r.h.s.\ in
the above ``mirror''  equality of cohomology groups.

\sip

Thus, for any Calabi-Yau manifold $X$ there are two
semi-infinite variations of Hodge structure, ${\mathsf
VHS}^B({X})$ and
${\mathsf VHS}^A({X})$, and two families of solutions,
 $\{\Phi^W_{\mathrm compl}(X)\}$ and
$\{\Phi_{\mathrm sympl}^{W'}(X)\}$,  to WDVV equations. In the
idealized situation when $X$ and $\hat{X}$ are dual torus
fibrations over the same Monge-Amp\`ere manifold \cite{KS,L,SYZ},
one has
$$
{\mathsf VHS}^A({X}) = {\mathsf VHS}^B(\hat{X}), \ \
\ \ {\mathsf VHS}^B({X}) = {\mathsf VHS}^A(\hat{X}),
$$ implying
$$
\Phi^W_{\mathrm compl}(X)=\Phi_{\mathrm
sympl}^{\hat{W}}(\hat{X}), \ \ \ \ \Phi^{W'}_{\mathrm sympl}(X)=
\Phi_{\mathrm compl}^{\hat{W}'}(X)
$$
for appropriately related
filtrations $(W,\hat{W})$ and $(W',\hat{W}')$.
To extend these equalities to an arbitrary  mirror pair $X$ and
$\hat{X}$, one has to find a conceptual way of  incorporating
instanton corrections (which
vanish for dual torus fibrations) into the definition of ${\mathsf
VHS}^A$.

\sip

As a purely algebraic exercise, we show in this paper that
semi-infinite variations of Hodge structure and the associated
construction of solutions of WDVV equations make sense for any
differential Gerstenhaber-Batalin-Vilkovisky (dGBV) algebra satisfying
Manin's axioms \cite{Ma}.

\sip

The paper is organized as follows. Section 2 gives an outline of
deformation theory and introduces the basic algebraic input.
In Sections 3 and 4 we construct ${\mathsf VHS}^A(\hat{X})$ and $\{\Phi^{W'}_{\mathrm
sympl}(X)\}$. In Section 5 we establish isomorphisms,
${\mathsf VHS}^A({X}) = {\mathsf VHS}^B(\hat{X})$ and
${\mathsf VHS}^B({X}) = {\mathsf VHS}^A(\hat{X})$,
for dual torus fibrations.  (Semi-infinite) Variations of a dGBV theme are
collected in Section 6.

\bip

\bip

\begin{center}
{\bf \S 2. An outline of deformation theory}
\end{center}

\bip

{\bf 2.1. Sign conventions.} In the deformation theory context, it
is more suitable to work with the {\em odd}\, version of the usual
notion of  differential   Lie superalgebra. By
definition, this is a $\Z_2$-graded vector space, $\fg=\fg_{\tlo}\oplus
\fg_{\tln}$, equipped with two odd linear maps
$$
d: \fg \rar \fg, \ \ \ \ \ \mbox{and} \ \ \ \ \ [\ \bullet \ ]: \fg\ot \fg \rar \fg,
$$
such that
\Bi
\item[(a)] $d^2=0$,\\
\item[(b)] $[a\bullet b]=-(-1)^{(\tl{a}+1)(\tl{b}+1)}[b \bullet a]$,\\
\item[(c)] $d[a\bullet b] = [da\bullet b] - (-1)^{\tl{a}}[a\bullet d b]$,\\
\item[(d)] $
           [a\bullet [b\bullet c]]= [[a\bullet b]\bullet c] +
           (-1)^{(\tl{a}+1)(\tl{b}+1)} [b\bullet [a\bullet c]$,
\Ei
for all $a,b,c\in \fg_{\tlo}\cup \fg_{\tln}$.

\bip

In many important examples, the $\Z_2$-grading in $\fg$ comes from an
underlying $\Z$-grading, i.e.
$\fg=\oplus_{i\in \Z}\fg^i$, $\fg_{\tlo}= \oplus_{i\, \mathrm even}\fg^i$,
$\fg_{\tln}= \oplus_{i\, \mathrm odd}\fg^i$,
and the basic operations satisfy
$d\fg^i\subset \fg^{i+1}$, $[\fg^i\bullet \fg^j]\subset \fg^{i+j-1}$.

\bip

Clearly, the parity change functor transforms this structure into the usual
structure of differential Lie superalgebra on the vector superspace $\Pi \fg$.
Thus the above notion is nothing new except slightly unusual sign conventions.

\bip

\sip

{\bf 2.2. Deformation theory.} Here is a latest guide for
constructing a versal  moduli space, $\cal M$,  of deformations of a
given mathematical structure:
\Bi
\item[Step 1:] Associate to the mathematical structure a
``controlling" differential graded (dg, for short)
Lie algebra $(\fg=\bigoplus_{i\in \Z} \fg^i, d, [\, \bullet\, ])$.\\

\item[Step 2:] Find a mini-versal {\em smooth}\, formal
 pointed dg-manifold $(M, \eth, *)$ (i.e.\  a triple
consisting of a smooth formal $ \Z$-graded manifold $M$, a point
$*\in M$, and an odd vector field $\eth$ on $M$ such that
$[\eth,\eth]=0$) which represents\footnote{In the sense that
$\Df_{\fg}(A,d_A)\simeq
 \mbox{Mor}((A,d_A,0)^{op},(M,\eth,*))$,
where $(A,d_A,0)^{op}$ is the representative of $(A,d_A)$ in
the opposite category ${\cal A}^{op}$.}
 the
deformation functor,
$$
\Ba{rccc}
{\Df}_{\fg}: & \left\{\Ba{c} \mbox{the category, $\cal A$,
   of }\\
                    \mbox{dg Artin algebras}
                    \Ea \right\}&
\lon & \left\{\mbox{the category of sets}\right\}\vspace{3mm}\\
& (A,d_A) & \lon & {\Df}^{\, *}_{\fg}(A,d_A),
\Ea
$$
$$
{\Df}_{\fg}(A,d_A):= \frac{
\left\{\Gamma\in (\fg\ot m_A)^2 \mid d\Gamma
+ d_A\Gamma + \frac{1}{2}[\Gamma\bullet \Gamma]=0\right\} }
{\exp{(\fg\ot m_{\cB})_{\tln}}}.
$$
Here $m_A$ stands for the maximal ideal in $A$, and
the quotient is taken with respect to the
following action of the
gauge group,
$$
\Gamma \rar \Gamma^g = e^{\ad_g}\Gamma -
\frac{e^{ad_g}-1}{\ad_g}(d+d_A)g,
 \ \ \ \ \forall g\in (\fg\ot m_A)^1.
$$

\sip

It is proven in \cite{Me2} that such a smooth  dg-manifold $(M,\eth,*)$ always
exists. Moreover, $(M,*)$ can be identified with a neighborhood of zero
in the cohomology superspace, $H(\fg)=\ker d/ \Img d$, so that the main job in Step 2
 is to find the vector field $\eth$.

\item[Step 3:]  Try to make sense to the quotient
space, $\cM$, of the subspace ${\mathsf zeros }(\eth) \subset M$ with respect to the
foliation governed by the integrable distribution $\Img \eth:=[\eth, TM]$,
$TM$ being the tangent
sheaf. In the analytic category, the Kuranishi
technique (that is, a cohomological splitting of $\fg$ induced by a suitably
chosen norm) should do the job.
This $\cM$ is a desired mini-versal moduli space of deformations.

\Ei

\sip

The functor $\Df_{\fg}$ is called {\em non-obstructed}\, if the vector field
 $\eth$ vanishes. In this case ${\cal M}=M$.

\bip
\sip

{\bf 2.3. Remarks.}  The classical deformation functor,
$\Def_{\fg}$,
as defined in the works of Deligne, Goldman, Kontsevich, Millson and others
(see  \cite{GM,Kon} and references therein),
 is just a restriction of $\Df$ to the subcategory, $\{(A, d_A=0)\}$, of
 usual (non differential) Artin superalgebras. Here is an evidence
 in support of the Cyrillic   version:
\Bi
\item[(i)] The modification $\Def_{\fg}\rar \Df_{\fg}$ does {\em
not}\, break the Main Theorem of Deformation Theory:
if dg Lie algebras $\fg_1$ and $\fg_2$ are quasi-isomorphic,
then $\Df_{\fg_1}\simeq \Df_{\fg_2}$.
\item[(ii)] Contrary to $\Def$, the functor $\Df$ is {\em
always}\, representable by  {\em smooth}\, geometric data $(M,\eth,*)$.
\item[(iii)] If $\Def_{\fg}$ admits a mini-versal
(usually, singular) moduli
space $\cM$, then the latter can be reconstructed from
$(M,\eth,*)$  as in Step 3.
In particular,  $\Def_{\fg}$ is non-obstructed in the usual sense
if and only if $\Df_{\fg}$ is non-obstructed.
\Ei

\sip

Rather than working with a singular moduli space ${\cal M}$,
it is more convenient to work with its smooth resolution
$(M,\eth,*)$.
This is the main, purely technical, advantage of $\Df$
over $\Def$.

\sip

Note  that the tangent space to the
functor $\Df_{\fg}$ is the full cohomology group $H(\fg)=\oplus_i H^i(\fg)$
rather than its subgroup $H^2(\fg)$ as in many classical
deformation problems (see examples below).
 This
is because we allowed the solutions, $\Gamma$, of Maurer-Cartan
equations to lie in $(\fg\ot m_A)^2$ rather than in $\fg^2\ot
m_A$. In this sense the  moduli space ${\cal M}$ (or its resolution
$(M,\eth,*)$)
describes {\em extended}\, deformations of the mathematical
structure under consideration. The classical moduli space, ${\cal
M}_{cl}$,
is a proper subspace of ${\cal M}$.

\bip

\sip

{\bf 2.4. Example (Deformations of complex structures).}
The dg Lie algebra controlling deformations of a given complex
structure on a $2n$-dimensional manifold $X$ is given by
$$
\fg=\left(\bigoplus_{i=0}^{2n} \fg^i,\   \fg^i= \bigoplus_{p+q=i}
\Gamma(M, \wedge^p T_X\ot \Omega^{0,q}_{X}),\ [\ \bullet\ ], \bp
\right)
$$
where $T_X$ stands for the sheaf of holomorphic vector fields,
$\Omega^{s,q}_X$ for the sheaf of smooth differential forms of
type $(s,q)$, and $[\ \bullet\ ]= {\mathsf Schouten\ brackets}\ot
{\mathsf wedge\ product}$.

\sip

In general, the deformation theory is obstructed. However,
if $X$ is a Calabi-Yau manifold, then $\Def_{\fg}$ (or
$\Df_{\fg}$)
is non-obstructed \cite{BK}, and the associated mini-versal moduli space
$\cM$ is isomorphic to an open neighbourhood of zero in
$H(\fg)=H^*(X,\wedge^* T_X)$. The embedding $\cM_{cl}\subset \cM$
corresponds to the inclusion $H^1(X,T_X)\subset H^*(X,\wedge^*
T_X)$.

\sip
\bip

{\bf 2.5. Example (deformations of Poisson and symplectic
structures).} The dg Lie algebra controlling deformations of a
given Poisson structure, $\nu_0\in \Gamma(X, \wedge^2 T_{\R})$, on
a real smooth manifold $X$ is given by
$$
\left(\oplus_{i=0}^{\dim X} \Gamma(X, \wedge^i T^*_{\R}),
[\ \bullet\ ]={\mathsf Schouten\ brackets}, d=[\nu_0\bullet\ldots ]
\right),
$$
where $T_{\R}$ stands for the sheaf of real tangent vectors.

\sip

If $\nu_0$ is non-degenerate, that is, $\nu_0=\omega^{-1}$
for some symplectic form $\omega$ on $X$, then the natural
``lowering of indices map'' $\omega^{\wedge i}:
\wedge^i T_{\R} \rar \wedge^i T_{\R}^*$ sends $[\nu_0\bullet\ldots ]$
into the usual de Rham differential. The image of the Schouten brackets
under this isomorphism we denote
by $[\ \bullet\ ]_{\omega}$. In this way we make the de Rham complex of
X
into a dg Lie algebra,
$$
\fg= \left(\bigoplus_{i=0}^{\dim X} \Gamma(X, \wedge^i T^*_{\R}),\
[\ \bullet\ ]_{\omega}, \
d={\mathsf de\ Rham\ differential}\right),
$$
which controls the extended deformations of the symplectic
structure $\omega$. More explicitly,
\[
 [\ka_1\bullet \ka_2]_{\omega}:= (-1)^{\tlk_1}[i_{\omega^{-1}}, d]
(\ka_1\wedge \ka_2) -
(-1)^{\tlk_1}\left([i_{\omega^{-1}}, d]\ka_1\right)\wedge \ka_2 -
\ka_1\wedge[i_{\omega^{-1}}, d]\ka_2, \ \ \ \forall \ka_1,\ka_2\in
\wedge^*T_{\R},
\]
with  $i_{\omega^{-1}}:\wedge^*T_{\R}\rar \wedge^{*-2}T_{\R} $ being
the natural contraction with the 2-vector $\omega^{-1}=\nu_0$.

\sip

It is proven in \cite{Me1}  that the associated deformation functor
is non-obstructed provided the symplectic manifold $(X,\omega)$ is
Lefschetz (in particular, K\"ahler), that is, that the cup
product
\[
[\omega^k]: H^{n-k}(X,\R) \lon H^{n+k}(X,\R)
\]
is an isomorphism for any $k\leq n=:1/2\dim X$. In this case the associated
mini-versal moduli space, $\cM$, of extended symplectic structures is smooth,
and  locally isomorphic to $H^*(X,\R)$.
The embedding $\cM_{cl}\subset \cM$
corresponds to the inclusion $H^2(X,\R)\subset H^*(X,\R)$.

\sip
\bip

{\bf 2.6. Example (extended deformations of K\"ahler
structures)}. Assume $X$ is a complex manifold, and
consider the Lie algebra,
$$
\fg'=\left(\bigoplus_{i=0}^{2n} \fg^i,\   \fg^i= \bigoplus_{p+q=i}
\Gamma(X, \wedge^p \overline{T}_X \ot \wedge^q T_X),\ [\ \bullet\ ]
={\mathsf Schouten\ brackets}\ot {\mathsf wedge\ product}
\right)
$$
If $\omega$ is a Kahler form
on $X$ and $\omega^{-1}\in \Gamma(X, \overline{T}_X \ot T_X) $
is its inverse, then $[\omega^{-1}\bullet\ldots ]$ defines a differential
on $\fg'$. Clearly, the resulting dg Lie algebra controls
deformations of the K\"ahler structure. It has a more convinient
embodiment.

\sip

The K\"ahler form $\omega$ induces the ``lowering of indices''
isomorphism,
$$
\Gamma(X,\wedge^p \bar{T}_X \ot \wedge^q T_X)
\lon \Gamma(X, \Omega^{p,q}_{X}),
$$
which sends $[\omega^{-1}\bullet\ldots]$ into $\bp$, and hence
makes the Doulbeout complex into a dg Lie algebra,
$$
\fg=\left(\bigoplus_{i=0}^{2n} \fg^i,\   \fg^i= \bigoplus_{p+q=i}
\Gamma(X, \Omega^{p,q}_{X}),\ [\ \bullet\ ]_{\omega},
\
\bar{\p}\right).
$$
The odd brackets can be written explicitly  as in 2.5 but with
$[i_{\omega^{-1}},d]$ replaced by $[i_{\omega^{-1}}, \bar{\p}]$.
The associated deformation functor is
non-obstructed, and the associated versal moduli space, $\cM$,
 of extended
K\"ahler forms is locally
isomorphic, as a formal pointed supermanifold, to
 $H^*(X,\C)=H^*(X,\Omega^*_Y)$, where $\Omega^*_X$
stands for
the sheaf
of holomorphic differential forms.
We often call $(\cM,*)$ the {\em extended K\"ahler cone}\,\footnote{This terminology
could be misleading as we do {\em not}\,  fix a particular isomorphism
$(M,*)\simeq (H^*(X,\C),0)$. Eventually, however, we will employ a
distinguished family of such  isomorphism  parameterized by isotropic
filtrations of $H^*(X,\wedge^*T_X)$ which are complementary to the
standard Hodge one.} of $X$.
The embedding $\cM_{cl}\subset \cM$
corresponds to the inclusion $H^1(X,\Omega^1)\subset H^*(X,\Omega^*_X)$.

\bip

\bip


\begin{center}
{\bf \S 3. A local system on the extended K\"ahler cone}
\end{center}

\bip

{\bf 3.1. From $\fg$-modules to vector bundles.} Let $(\fm,
\bullet, d)$ be a dg module over the dg Lie algebra $(\fg, [\
\bullet\ ], d)$, that is, a $\Z$-graded vector space $\fm$ together with
two odd linear maps, $d:\fm \rar \fm$ and $\bullet : \fm\ot \fg
\rar \fm$ such that
\Bi
\item[(i)] $d^2=0$,
\item[(ii)] $d (\ka\bullet a) = (d \ka)\bullet a -
(-1)^{\tl{\ka}}\ka\bullet da$,
\item[(iii)] $\ka_1\bullet \ka_2 \bullet a -
(-1)^{(\tlk_1+1)(\tlk_2+1)} \ka_2\bullet \ka_1 \bullet a
=[\ka_1\bullet \ka_2] \bullet a$
\Ei
for all $\ka_1,\ka_2\in \fg$ and $a\in \fm$.

\sip

Let $(M,*,\eth)$ be the mini-versal dg moduli
space associated with the functor $\Df_{\fg}$.

\sip

The dg Lie algebra structure on the vector superspace $\fg$ can be
geometrically represented as an odd homological vector field,
$Q_{\fg}$, on $\fg$ viewed as a supermanifold. For any $\al \in
\fg^*$ (interpreted now as a function on the supermanifold $\fg$)
and any point $\ga\in \fg$, one has $$ Q_{\fg} \al\mid_{\ga} =
-(-1)^{\tl{\al}}\langle\al, d\ga + \frac{1}{2}[\ga\bullet
\ga]\rangle. $$

\sip

There always exists a map of pointed dg manifolds,
$$
\Gamma:
(M,*,\eth) \lon (\fg, 0, Q_{\fg})
$$
such that $d\Gamma: T_*M \rar
T_0\fg$ is a monomorphism. Moreover, this map is unique up to a
gauge transformation as in Step 2 of Sect.\ 2.2 (see, e.g., \cite{Me2}
for a proof).

\sip

Any such $\Gamma$ gives rise to a flat $\eth$-connection,
$$ \Ba{rccc}
D_{\eth}^{\Gamma}: & \f_{M} \ot \fm & \lon & \f_M\ot \fm\\
                 &  f a      & \lon & D_{\eth}^{\Gamma}(f a):=
                 \eth f a + (-1)^{\tl{f}} f (da +
                 \Gamma\bullet a),

\Ea
$$
 in the trivial vector bundle
$M\times \fm$. Though this connection is not gauge invariant,
the associated cohomology sheaf,
$$
E_{\fm}:= \f_M \ot_{\Ker \eth} \frac{\Ker D_{\eth}^{\Gamma}}{\Img
D_{\eth}^{\Gamma}},
$$
together with its natural flat
$\eth$-connection, $D_{\eth}$,  is well defined, i.e.\ the pair
$(E_{\fm},\eth)$ does not
depend on the choice of a particular map $\Gamma$.

\sip


In summary, we have the following

\bip

{\bf 3.1.1. Proposition.} {\em  Let $(\fg, [\ \bullet\ ], d)$ be a
dg Lie algebra
and let $(M,*,\eth)$ be the mini-versal dg moduli
space representing the functor $\Df_{\fg}$. Any dg\, $\fg$-module
$(\fm, \bullet, d)$ gives canonically rise
\Bi
\item[(i)] to an $\f_M$-module, $\pi: E_{\fm}\rar M$,
such that $\pi^{-1}(*)=H(\fm)$, the cohomology of the complex
$\fm$, and
\item[(ii)] to a flat $\eth$-connection $D_{\eth}: E_{\fm}\rar
E_{\fm}$.
\Ei }

\bip

It is clear that the pair $(E_{\fm}, D_{\eth})$ on $M$ gives rise
to a well defined $\f_{\cM}$-module on the derived moduli space
$\cM={\mathsf Zeros}(\eth)/ \Img \eth$ (see \cite{Ba}).

\bip

{\bf 3.1.2. Remark.} The above Proposition can be strengthened
as follows: the  derived category of dg
modules over a given dg Lie algebra $\fg$ is equivalent to
the purely geometric category of vector bundles
over $(M,\eth,*)$ equipped with flat $\eth$-connections. We omit the
proof.

\sip

In fact the base, $(M,\eth,*)$, can
itself be identified with  $F(\fg)$, where $F$ is the functor
$$
\left\{ \Ba{c} \mathrm the\ category\ of\\ \mathrm dg\ Lie\ algebras
         \Ea\right\}
         \stackrel{F}{\lon}
\left\{ \Ba{c} \mathrm the\ derived\ category\ of\\ \mathrm dg\ Lie\ algebras
         \Ea\right\}.
$$
This gives another meaning
to the dg moduli space $(M,\eth,*)$ which we first encountered in the
context of the $\Df$ormation theory.

\bip

\sip

{\bf 3.2. Flat structure.} Let $(\fm, d, \bullet)$ be a dg
$\fg$-module and assume that there is an even linear map,
$$
\Ba{rccc}
\circ: & \fg\ot \fm & \lon & \fm\\
        &  \ka \ot a & \lon & \ka\circ a
\Ea
$$
such that
$$
\ka_1\circ \ka_2 \circ  a = (-1)^{\tlk_1\tlk_2}\ka_2\circ \ka_1 \circ a
,
$$
$$
d (\ka\circ a)
=
  (d \ka) \circ a + (-1)^{\tlk}\ka
 \circ d a + (-1)^{\tl{\ka}} \ka \bullet a,
$$
and
$$
\ka_1\circ \ka_2 \bullet a - (-1)^{\tlk_1\tlk_2 + \tlk_1} \ka_2\bullet
\ka_1\circ a = -(-1)^{\tlk_2}[\ka_1\bullet \ka_2]\circ a.
$$
for any $\ka,\ka_1,\ka_2\in \fg$ and $a\in \fm$.

\sip

Choosing a  map $\Ga: (M,*,\p)\rar (\fg,0,Q_{\fg})$ as
in Sect.\ 3.1, we may define  a flat connection,
$$
\Ba{rccc}
\nabla^{\Gamma}:& TM \ot (\f_M \ot\fm) & \lon & \f_M\ot \fm \\
     &    v \ot (f a) & \lon & \nabla^{\Gamma}_v(fa):= (v f) a + (-1)^{\tl{v}\tl{f}}
         f\left((v\Gamma)\circ a\right),
\Ea
$$
in the trivial vector bundle $M\times \fm$. It is not hard to
check that
$$
e^{-\Gamma\circ}(d+\eth) e^{\Gamma\circ} = D^{\Gamma}_{\eth},
$$
and
$$
e^{-\Gamma\circ}(v) e^{\Gamma\circ} = \nabla^{\Gamma}_{v},
$$
for any vector field  $v$ on $M$. Hence,
$$
[\nabla_v^{\Gamma}, D_{\eth}^{\Gamma}] = \nabla_{[v,\eth]}^{\Gamma}
$$
implying the following


\bip

{\bf 3.2.1. Proposition} \cite{Ba}. {\em  Let $(\fg, [\ \bullet\ ], d)$ be a
dg Lie algebra, $(M,*,\eth)$ the associated mini-versal dg moduli
space, and $\cM={\mathsf Zeros}(\eth)/ {\mathsf Im}(\eth)$ the associted derived moduli space.
Any dg $\fg$-module
$(\fm, \bullet, \circ, d)$ as above gives canonically rise
to a pair, $(E_{\fm}, \nabla)$, where $E_{\fm}$ is a vector bundle
on $\cM$ with typical fibre $H(\fm)$,
and $\nabla$ is a flat connection.
}

\bip

{\bf 3.2.2. Remark.} The flat connection $\nabla$ identifies the
linear space of horizontal sections of $\pi: E_{\fm}\rar \cM$ with the
fibre $\pi^{-1}(*)$. Explicitly, the identification goes as
follows
$$
\Ba{ccc}
H(\fm) & \stackrel{1:1}{\lon} & {\mathrm Space\ of\ horizontal\
sections}\\
a & \lon & e^{-\Gamma\circ} a.
\Ea
$$

\sip

In particular, the parallel transport establishes a canonical
isomorphism of the fibres $\pi^{-1}(t)\simeq \Ker (d+  \Gamma\bullet)/\Img
(d+ \Gamma\bullet)$ with the fibre, $\pi^{-1}(*)\simeq \Ker d/\Img
d$, over the base point.

\bip
\sip

{\bf 3.3. Example (deformations of complex structures)}.
One of the key observations in \cite{Ba} is that, for any
compact complex manifold $X$,
the pair consisting of the dg Lie algebra,
$$
\fg=\left(\bigoplus_{i=0}^{2n} \fg^i,\   \fg^i= \bigoplus_{p+q=i}
\Gamma(M, \wedge^p T_X\ot \Omega^{0,q}_{X}),\ [\ \bullet\ ]_{\mathrm Sch},
\bp
\right)
$$
and the $\fg$-module,
$$
\fm=\left(\bigoplus_{i=0}^{2n} \fm^i,\   \fg^i= \bigoplus_{p+q=i}
\Gamma(X, \Omega^{p,q}_{X}),\ \bullet\ ,\ \circ \ ,
d=\p+\bar{\p}\right),
$$
with $\ka\bullet a:= [i_{\ka}, \p]$ and $\ka\circ a:= i_{\ka}a$,
do satisfy the conditions of Proposition~3.2.1. Thus the
extended moduli
space of complex structures (which, for Calabi-Yau $X$,  is smooth
and  locally isomorphic to
$H^*(X,\wedge^*T_X)$)
comes equipped with a flat vector bundle whose typical fibre
is the de Rham cohomology $H^*(X,\C)$.

\sip

\bip

{\bf 3.4. Example (deformations of K\"ahler forms).}
In this subsection we shall present one more example to which
Proposition~3.2.1 is applicable. Curiously it inverses
the roles of $\fg$ and $\fm$ in Barannikov's example, and
hence gives rise to a local
system whose base is locally isomorphic to the de Rham
cohomology $H^*(X,\C)$ and whose typical fibre is
$H^*(X,\wedge^* T_X)$.

\sip

We assume from now on that $(X,\omega)$ is
a K\"ahler manifold, and denote by $\fg$ the dg Lie algebra
controlling extended deformations of the K\"ahler structure (see
Sect.\ 2.6). We omit from now on the subscript $\omega$ in the Lie
brackets notation.

\sip

\bip

{\bf 3.4.1. Auxiliary operators.}
We want to study the following morphisms of sheaves:
\begin{itemize}
\item[(i)] The inverse K\"ahler form $\omega^{-1}$ induces a natural ``raising of
indices'' map
\[
\sharp: \Omega^{0,q}_X \lon T_X\ot\Omega^{0,q-1},
\]
which, combined with the antisymmetrisation, extends to the map
\[
\sharp: \wedge^* T_X\ot \Omega^{0,*}_X \lon \wedge^{*+1}
T_X\ot\Omega^{0,*-1}_X.
\]
\item[(ii)]
For any $\kappa\in \Gamma(X, \Omega^{s,t}_X)$ there is a natural
map,
\[
\Ba{rccc}
 i_{\kappa}: & \wedge^* T_X\ot \Omega^{0,*}_X & \lon
  \wedge^{*-s} T_X\ot \Omega^{0,*+t}_X
\Ea
\]
\end{itemize}
which is a combination of contraction and wedge product.



\sip
\bip

{\bf 3.4.2. Lemma.} {\em For any K\"ahler manifold $X$, one has
\begin{itemize}
\item[(a)] The commutator,
$$
Q:=[ \sharp\, , \bp ]: \wedge^* T_X\ot \Omega^{0,*}_X \lon
       \wedge^{*+1} T_X\ot\Omega^{0,*},
$$
is a differential, i.e. $Q^2=0$.
\item[(b)] $[Q, \bp]=0$.
\item[(c)] $[Q, \sharp\,]=0$.
\item[(d)] $[i_{\ka_1}, [Q, i_{\ka_2}]]=-i_{[\ka_1\bullet \ka_2]}$
\ \ \ $\forall \ \ka_1,\ka_2 \in \fg$.
\item[(e)] $[i_{\ka_1\wedge\ka_2}, Q]= i_{\ka_1}[i_{\ka_2},Q] +
(-1)^{\tlk_1\tlk_2}i_{\ka_2}[i_{\ka_1}, Q]$
\ \ \ $\forall \ \ka_1,\ka_2 \in \fg$.
\end{itemize}
}

\sip

{\bf A comment on the proof.} One can identify the sheaf
$\Omega^{*,*}_{X}$ with the structure sheaf on the supermanifold
$\Pi T_{\C}$, where $\Pi$ is the parity change functor and
$T_{\C}=T_{\R}\ot \C$.  Then,
in a natural local coordinate system
$(z^a, \psi^a:= dz^a, \bar{\psi}^{a}:= d \bar{z}^a)$,
one may explicitly represent the (odd) Lie brackets  as  Poisson ones,
$$
[\ka_1 \bullet \ka_2] = (-1)^{\tlk_1} \omega^{a\bar{b}}\frac{\p \ka_1}
{\p \psi^a}\frac{\p \ka_2}{\p \bar{z}^b}  -
\bar{\psi}^{{c}} \frac{\p \omega^{a\bar{b}}}{\p \bar{z}^c}
\frac{\p \ka_1}{\p \psi^a} \frac{\p \ka_2}
{\p \bar{\psi}^{{b}}}  - (-1)^{(\tlk_1+1)(\tlk_2+1)}\left( \ka_1
\leftrightarrow
\ka_2\right)
$$
where the summation over repeated indices is assumed, and $\omega^{a\bar{b}}$
stand for the coordinate components of the inverse K\"ahler form $\omega^{-1}$.

\sip

Analogously, one can identify the sheaf $\wedge^* T_X\ot \Omega^{0,*}_X$
 with a  sheaf of functions on the total superspace
 of the bundle $\Pi T_X \oplus \Pi \overline{T}_X^*$.
 In a natural local coordinate chart,
 $(z^a, \psi_a:= \Pi \p/\p z^a, \bar{\psi}^a:= d\bar{z}^a)$,
 one may explicitly represent the basic operators as follows
 \begin{eqnarray*}
 \bp &=& \bar{\psi}^a   \frac{\p}{\p \bar{z}^a},\\
\sharp&=& \omega^{a\bar{b}}\psi_a\frac{\p}{\p \bar{\psi}^b} ,\\
Q &=& \omega^{a\bar{b}}\psi_a \frac{\p}{\p \bar{z}^b}
- \bar{\psi}^{{c}} \frac{\p \omega^{a\bar{b}}}{\p \bar{z}^c} \psi_a
\frac{\p}{\p \bar{\psi}^b}, \\
i_{\ka} &=& (-1)^{pq} \ka_{a_1\ldots a_p \bar{b}_1 \ldots \bar{b}_q}
\bar{\psi}^{b_1}\ldots \bar{\psi}^{b_q} \frac{\p^p}{
\p \psi_{a_1} \cdots \p\psi_{a_p}}, \ \ \  \forall \ka\in \Omega^{p,q}_M.
 \end{eqnarray*}

\sip

The main  technical advantage of this point of view is that
\begin{itemize}
\item[-] it is enough to check
all the claims (a)-(e) at one arbitrary point $*\in X$,
\item[-] due to the K\"ahler
condition on $\omega$, one can always
choose the coordinates at $*$ in such a way that
 $$
 \frac{\p \omega^{a\bar{b}}}{\p \bar{z}^c}\mid_{*}=0.
 $$
\end{itemize}

 With this observation all the above  expressions can be dramatically
 simplified making the claims either transparent or
  requiring a minimal calculation.
 \hfill $\Box$

\sip

\bip

{\bf 3.4.3. Proposition.}  {\em The (odd) linear map,
\[
\Ba{rccc}
\bullet: &  \fg \ot \fm  &  \lon \fm \\
                   &   \ka \ot a & \lon -[i_{\ka}, Q]a
\Ea
\]
makes the dg vector space,
\[
\fm:=\left(\bigoplus_{i=0}^{2n} \fm^i,\ \fm^i= \bigoplus_{p+q=i}
\Gamma(X, \Lambda^p T_X\ot \Omega^{0,q}_X),\
d=\bar{\p}+Q\right).
\]
into a dg module over the dg Lie algebra $\fg$.}

\sip

\Proof  For any  $\kappa_1, \kappa_2\in \fg$ and any $a\in \fm$, we have,
by Lemma~2.1(d),
\Beqrn
 \ka_1\bullet \ka_2 \bullet a + (-1)^{\tlk_1\tlk_2 + \tlk_1+\tlk_2}
 \ka_2\bullet \ka_1 \bullet a &=& [i_{\ka_1},Q]  [i_{\ka_2},Q]a
+ (-1)^{\tlk_1\tlk_2 + \tlk_1+\tlk_2} [i_{\ka_1},Q]  [i_{\ka_2},Q]a \\
&=& \left[\left[i_{\ka_1}, \left[Q,i_{\ka_2}\right]\right],Q\right]a\\
&=&-\left[i_{[\ka_1\bullet \ka_2]},Q\right]a\\
&=& [\ka_1\bullet\ka_2]\bullet a.
\Eeqrn
which means that $(\fm, \bullet)$ is a $\fg$-module. Consistency
 of $\bullet$ with the differentials is also an easy check:
\Beqrn
d(\ka\bullet a) &=& (\bp +Q)[i_{\ka},Q]a\\
&=&\bp i_{\ka} Qa -(-1)^{\tlk}\bp Qi_{\ka}a + Qi_{\ka}Qa\\
&=& [\bp, i_{\ka}]Qa - (-1)^{\tlk}i_{\ka}Q\bp a
+ (-1)^{\tlk}Q[\bp, i_{\ka}]a + Qi_{\ka}\bp a  - (-1)^{\tlk}[i_{\tlk},Q]Qa\\
&=& [i_{\bp \ka},Q]a - (-1)^{\tlk}[i_{\tlk},Q]Qa - [i_{\ka},Q]\bp a\\
&=&(\bp \ka)\bullet a - (-1)^{\tlk}\ka \bullet da. \hspace{9cm}  \Box
\Eeqrn



\bip
\sip
{\bf 3.4.4. Lemma.}
{\em An (even) linear
map
\[
\Ba{rccc}
 \circ: &    \fg \ot \fm  &  \lon & \fm \\
                   &   \ka \ot a &\lon &  i_{\ka}a
\Ea
\]
satisfies,

$$
d (\ka\circ a)=
  (\bp \ka) \circ a + (-1)^{\tlk}\ka
 \circ d a + (-1)^{\tl{\ka}} \ka \bullet a,
$$
and
$$
\ka_1\circ \ka_2 \bullet a - (-1)^{\tlk_1\tlk_2 + \tlk_1} \ka_2\bullet
\ka_1\circ a = -(-1)^{\tlk_2}[\ka_1\bullet \ka_2]\circ a.
$$
for any $\ka,\ka_1,\ka_2\in \fg$ and $a\in \fm$.

}

\bip

\Proof We have
\begin{eqnarray*}
d (\ka\circ a) &=& (\bp + Q)i_{\ka} a \\
 &=& [\bp,i_{\ka}]  a + (-1)^{\tl{\ka}}i_{\ka}\bp  a +
 [Q,i_{\ka}]a  + (-1)^{\tlk} i_{\ka}Qa\\
&=& i_{\bp \ka} a +(-1)^{\tl{\ka}}i_{\ka}(\bp+Q)  a - (-1)^{\tlk} [Q,i_{\ka}]a\\
&=&  (\bp \ka) \circ a + (-1)^{\tlk}\ka
    \circ d a + (-1)^{\tl{\ka}} \ka \bullet a,
\end{eqnarray*}
and, using Lemmma 2.1(d),
\Beqrn
\ka_1\circ \ka_2 \bullet a - (-1)^{\tlk_1\tlk_2 + \tlk_1} \ka_2\bullet
\ka_1\circ a &=& -i_{\ka_1}[i_{\ka_2},Q]a + (-1)^{\tlk_1\tlk_2 + \tlk_1}
[i_{\ka_2},Q]i_{\ka_1}a \\
&=& (-1)^{\tlk_2}[i_{\ka_1},[Q,i_{\ka_2}]]a\\
&=&  -(-1)^{\tlk_2}i_{[\ka_1\bullet \ka_2]}a\\
&=& -(-1)^{\tlk_2}[\ka_1\bullet \ka_2]\circ a.
\Eeqrn
\hfill $\Box$

\bip

In conclusion, the pair consisting of the dg Lie algebra,
$$
\fg=\left(\bigoplus_{i=0}^{2n} \fg^i,\   \fg^i= \bigoplus_{p+q=i}
\Gamma(X, \Omega^{p,q}_{X}),\ [\ \bullet\ ]_{\omega},
\
\bar{\p}\right),
$$
and the dg $\fg$-module,
\[
\fm:=\left(\bigoplus_{i=0}^{2n} \fm^i,\ \fm^i= \bigoplus_{p+q=i}
\Gamma(X, \Lambda^p T_X\ot \Omega^{0,q}_X),\ \bullet\ , \ \circ
\ , d=\bar{\p}+Q\right),
\]
satisfy the conditions of Proposition~3.2.1, and hence gives rise
to a flat vector bundle, $(E_{\fm}, \nabla)$, over the moduli space
of extended K\"ahler structures, $\cM=M\simeq H^*(X,\C)$, with
typical fibre $H^*(X,\wedge^*T_X)$.

\bip

{\bf 3.5. Remark.} In Sect.\ 5 we shall give more examples of
pairs $(\fg, \fm)$ to which Proposition 3.2.1 is applicable ---
one for each differential Gerstenhaber-Batalin-Vilkoviski (dGBV,
for short) algebra. This will produce, in particular, one more
local system, $(E'_{\fm}, \nabla)$, over the extended K\"ahler
cone of $X$, with typical fibre $H^*(X,\C)$. If $X$ is Ricci flat,
then $(E'_{\fm}, \nabla)$
 is isomorphic to the one constructed in subsect.\ 3.4,
but in general it is different.

\bip

 \bip

\pagebreak
\begin{center}
{\bf \S 4. \ \  Frobenius manifolds from semi-infinite
variations\\ of Hodge structure in $H^*(X,\wedge^*T_X)$}
\end{center}

\bip

\sip

{\bf 4.1. Semi-infinite $B$-variations of Hodge structure
in $H^*(X,\C)$.} A complex structure $J_0$
on  a compact  manifold
$X$ gives rise to the Hodge decomposition,
$\bigoplus_{i,j} H^i(X,\Omega^j)$,
of the de Rham cohomology group $H^*(X,\C)$. When the complex
structure is deformed,  the associated
Hodge filtration $F^{\geq r}_0$  gets  deformed
into another one, $F^{\geq r}_t$, $t\in \cM_{cl}$;
a remarkable fact is that this deformation
satisfies Griffiths transversality
condition with respect to the Gauss-Manin connection
on the bundle $E_{\fm}$ (see Examples 2.4 and 3.3)
restricted to $\cM_{cl}\subset
\cM$.

\sip

What happens to the Hodge filtration when one
moves from a ``classical''
point $J_t\in \cM_{cl}$ to a generic
point in the extended moduli space, $\cM$,
 of complex structures?
An answer to this question was given in \cite{Ba}
by a creative usage of Sato type Grassmanian,
$Gr_{\frac{\infty}{2}}$,
of semi-infinite subspaces in
$H^*(X,\C)[[\hbar, \hbar^{-1}]]$.
Moreover, for Calabi-Yau $X$, the resulting datum was used as an input
for producing a family of Frobenius manifold structures
on $\cM\simeq H^*(X,\wedge^* T_X)$ parametrized
by isotropic (with respect to the Poincare
metric on $H^*(X,\C)$) increasing filtrations which are
complementary to the Hodge one.

\sip

We present below a symplectic version
of the Barannikov's construction.
The parallelism is so strong that we can afford being
sketchy.

\sip
\bip

{\bf 4.2. Semi-infinite $A$-variations of Hodge structure
in $H^*(X,\wedge^* T_X)$.} Let $(X,\omega)$ be a compact $n$-dimensional
K\"ahler
manifold, and let $(\cM, *)\simeq (H^*(X,\C), 0)$ be the associated
moduli space of extended K\"ahler structures.
As before, we denote by $\fg$ the dg Lie algebra,
$$
\left(\bigoplus_{i=0}^{2n} \fg^i,\   \fg^i= \bigoplus_{p+q=i}
\Gamma(X, \Omega^{p,q}_{X}),\ [\ \bullet\ ]_{\omega},
\
\bar{\p}\right),
$$
which controls deformations of the K\"ahler structure.
However, instead of  the $\fg$-module $\fm$
defined in Sect.\ 3.4, we need its slight modification
involving a formal parameter $\hbar$,
\[
\fm:=\left(
\Gamma(X, \Lambda^* T_X\ot \Omega^{0,*}_X)[[\hbar,
\hbar^{-1}]],\ \bullet\ , \ \circ
\ , d=\bar{\p}+\hbar Q\right),
\]
with $\ka\bullet a:= -[i_{\ka},Q]a$ and $\ka\circ a
:= \frac{1}{\hbar} i_{\ka}a$.
All the formulae
and claims of Sect.\ 3.4 remain true; in particular,
the pair $(\fg,\fm)$ gives rise to a flat vector bundle,
$(E_{\hbar}, \nabla)$, over the moduli space
$(\cM,*)$. As
$$
\bp + \hbar Q= l_{\hbar}^{-1} \hbar^{\frac{1}{2}}(\p +Q) l_{\hbar},
$$
where $l_{\hbar}$ is a linear automorphism of
$\Gamma(X,\wedge^* T_X\ot \Omega^{0,*}_X )[[\hbar,\hbar^{-1}]]$
 given by\footnote{The associated linear automorphism
of $H^*(X,\wedge^* T_X)[[\hbar,\hbar^{-1}]]$ is denoted by the
same letter $l_{\hbar}$.}
$$
l_{\hbar}a:= \hbar^{\frac{n+q-p}{2}}a, \ \  \forall \ a\in \Gamma(X,\wedge^p
T_X\ot \Omega^{0,*}),
$$
the fibre of the bundle $E_{\hbar}$
over the base point
$*$ is isomorphic to $H^*(X,\wedge^* T_X)[[\hbar,\hbar^{-1}]]$.

\sip

The standard Hodge decreasing filtration on $H^*(X,\wedge^*T_X)$,
$$
\left(0=F^{\geq\frac{n+2}{2}}\subset F^{\geq\frac{n}{2}}\subset \ldots
\subset F^{\geq\frac{-n}{2}}\
, \ 0=F^{\geq\frac{n+1}{2}}\subset F^{\geq\frac{n-1}{2}}\subset
\ldots
\subset F^{\geq\frac{1-n}{2}} \right)
$$
is given by
$$
F^{\geq r}:= \bigoplus_{p-q\geq 2r \atop  \tl{p}=\tl{q} +\widetilde{2r}}
H^q(X,\wedge^p T_X), \ \ \ r\in \Z[\frac{1}{2}]
$$
(so that $F^{\geq \frac{n}{2}}=H^0(X,\wedge^n T_X)$,
$F^{\geq \frac{n-1}{2}}=H^0(X,\wedge^{n-1}T_X)\oplus H^1(X,\wedge^n
T_X)$,
etc.).
There is associated a linear subspace,
\Beqrn
L_0 &:=& \mbox{span}_{r\in \Z[\frac{1}{2}]}\, F^{\geq r}\hbar^{-r+ \frac{n}{2}}[[\hbar]]
\\  &=&\mbox{span}\, l_{\hbar} H^*(X,\wedge^* T_X)[[\hbar]]
\subset H^*(X,\wedge^* T_X)[[\hbar,\hbar^{-1}]].
\Eeqrn

\sip

Next one considers a
relative Grassmanian,
$Gr_{\frac{\infty}{2}}(E_{\hbar})$,
whose fibre over a generic point $t\in \cM$ is the Sato
 Grassmaninan of
semi-infinite subspaces in the fibre of   $\pi_{\hbar}: E_{\hbar}\rar \cM$ over $t$.
The latter has a canonical global section, $\Phi$,
which associates to any $t\in \cM$ the vector space, $\Phi(t)$, of all elements in
$\pi_{\hbar}^{-1}(t)$ which are (formal) analytic in $\hbar$.
 We can use the flat
connection $\nabla$ to compare   the values $\Phi(t)$ at different points
 $t\in \cM$ via the parallel transport to the base point.
In this  way we get a well defined map to the projective limit of Sato Grassmanians,
$$
\Ba{ccc}
 \cM & \lon & Gr_{\frac{\infty}{2}}(H^*(X,\wedge^* T_X)
[[\hbar,\hbar^{-1}]])(\f_*)\\
 t & \lon & L_t,
 \Ea
$$
where $L_t$ is the $\f_*$-submodule of
$H^*(X,\wedge^* T_X)[[\hbar,\hbar^{-1}]]\ot \f_*$ generated
by $l_{\hbar}P^{\nabla}_t \Phi(t)$, where $P^{\nabla}_t$ stands
for the parallel transport from $t$ to $*$. Taking into account
Remark 3.2.2, we may write $L_t= l_{\hbar}e^{\Gamma(t)\circ}
\Phi(t)\ot \f_*$. Note that $L_*$ is
precisely the submodule $L_0\ot \f_*$ corresponding to the
standard Hodge filtration on $H^*(X,\wedge^* T_X)$.

\bip


\sip

\bip

{\bf 4.3. Frobenius manifolds.} We assume from now on that $X$ is
 an $n$-dimensional Ricci flat K\"ahler manifold, i.e.\ a Calabi-Yau manifold. In
 this case the line bundle $\wedge^n T_X$ admits a global nowhere
 vanishing section which we denote by $\eta$. The associated
 nowhere vanhishing global holomorphic $n$-form is denoted by
 $\Omega$. Note that $\eta \in L_0$.

 \sip

 Let
$$
\left(0=W_{\leq\frac{-n}{2}}\subset W_{\leq\frac{-n+2}{2}}\subset \ldots
\subset W_{\leq\frac{n+2}{2}}\
, \ 0=W_{\leq\frac{1-n}{2}}\subset W_{\leq\frac{3-n}{2}}\subset
\ldots
\subset
W_{\leq\frac{n+1}{2}}
\right),
$$
 be an increasing
 filtration on $H^*(X, \wedge^* T_X)$ which is complementary
to the Hodge filtration in the sense that
$$
\bigoplus_{\tl{i}+\tl{j}=\widetilde{2r}} H^i(X, \wedge^i T_X) =
F^{\geq r}\oplus W_{\leq r}.
$$
There is associated a linear subspace,
$$
L_W := \mbox{span}_{r\in \Z[\frac{1}{2}]}\, W_{\leq r}\hbar^{-r+ \frac{n}{2}}[[\hbar^{-1}]]
\subset H^*(X,\wedge^* T_X)[[\hbar,\hbar^{-1}]].
$$
Moreover,
$$
H^*(X,\wedge^* T_X)[[\hbar,\hbar^{-1}]] = L_0 \oplus L_W
$$
and, for $t\in \cM$ ``sufficiently close'' to
the base point $*$, the intersection $L_t\cap (\eta + L_W)$
consists of a single element. Hence there is a well defined composition

$$
\Ba{rccccc}
\Psi^W: & \cM & \lon &  H^*(X,\wedge^* T_X)[[\hbar,\hbar^{-1}]]
& \stackrel{\lrcorner \Omega}{\lon} H^*(X, \Omega^{n-*})[[\hbar,\hbar^{-1}]]\\
        & t & \lon & L_t\cap (\eta +
        L_W) & \lon \left(L_t\cap (\eta +
        L_W)\right)\lrcorner \Omega
\Ea
$$
where $\lrcorner \Omega$ stands for the natural contraction
with the holomorphic volume form. Moreover its equivalence class,
$$
\hat{\Psi}^W:= (\Psi^W - \eta)\ \bmod
\hbar^{-1}L_W \lrcorner \Omega,
$$
 gives rise to a composition
$$
\hat{\Psi}^W: \cM \lon \left( L_W \bmod
\hbar^{-1}L_W\right)\lrcorner\Omega = \left(\bigoplus_r W_{\leq r}
/W_{\leq{r-1}}\right)\lrcorner\Omega = H^*(X,\C),
$$
which is
obviously a local diffeomorphism. If $\{\Delta_a\}$ is a vector
space basis in $H^*(X,\C)$ and $\{t^a\}$ the associated dual
basis, then the functions $t^a_W:=(\hat{\Psi}^W)^{-1}(t^a)$ define
a distinguished coordinate system, $\{t_W^a\}$, on
$\cM$\footnote{More precisely, the map $\hat{\Psi}^W$ defines a
distinguished flat structure, $\nabla$, on $\cM$.}. In this
coordinate system the map $\Psi^W$ satisfies the equations
(Proposition 6.5 in \cite{Ba}),
$$
\frac{\p^2 \Psi^W(t_W,
\hbar)}{\p t^a_W\p t^b_W} = \hbar^{-1}\sum_c A^c_{ab}(t_W)\frac{\p
\Psi^W(t_W, \hbar)}{\p t^c_W},
$$
and hence gives rise to a
pencil, $\p/\p t_W^a + \lambda A_{ab}^c(t_W)$, of flat connections
on the moduli space of extended K\"ahler structures. That is, the
{\em structure functions}\, $ A^c_{ab}(t_W)$ satisfy, $$ \frac{\p
A^c_{ab}}{\p t_W^d}= (-1)^{\tl{a}\tl{d}}\frac{\p  A^c_{db}}{\p
t_W^a}, \ \ \ \ \ \ \ \  \sum_c A_{ab}^c A_{cd}^e =
(-1)^{\tl{a}(\tl{b} + \tl{d})}
 \sum_c A_{bd}^c A_{ca}^e.
$$
Thus the product in the tangent sheaf ${\cal T}_{\cM}$ defined by
$$
\frac{\p }{\p t^a_W}\circ \frac{\p}{\p t^b_W}:= \sum_c A_{ab}^c
\frac{\p}{\p t^c_W}
$$
is associative.

\sip

The Poincare form on $H^*(X,\C)$ together with the given holomorphic
volume form $\Omega$ induce a non-degenerate paring on
$H^*(X,\wedge^*T_X)[[\hbar,\hbar^{-1}]]$,
$$
(v, w):= \int_X (v\lrcorner \Omega) \wedge (w\lrcorner \Omega).
$$
If we assume that the complementary filtration $W$ is isotropic in
the sense that
$$
 (W_{\leq r}, W_{\leq -r+1})=0,\ \ \ \ \forall\ r\in \Z[\frac{1}{2}],
$$
then,
$$
(v,w)\in \hbar^{n-2}[[\hbar^{-1}]], \ \ \ \ \ \forall\ v,w\in L_W.
$$
On the other hand,
$$
\int_X  \frac{\p\Psi^W(t_W, \hbar)}{\p t^a_W} \frac{\p\Psi^W(t_W, -\hbar)}{\p t^b_W}
\in \hbar^{n-2}[[\hbar]],
$$
so that the functions
$$
g_{ab}(t_W)
:= \hbar^{2-n}\int_X
\frac{\p\Psi^W(t_W, \hbar)}{\p t^a_W} \frac{\p\Psi^W(t_W, -\hbar)}{\p t^b_W}
$$
do not depend on $\hbar$, and hence define a non-degenerate
metric, $\langle \frac{\p }{\p t_W^a}, \frac{\p}{\p t_W^b}\rangle
$, on the formal pointed manifold $(\cM,*)$. Moreover,
differentiating the above equation with respect to $t^c_W$, we
immediately conclude that the metric $g_{ab}$ is constant in these
coordinates , and that the tensor,
$$
A_{abc}:= A_{ab}^dg_{dc},
$$
is totally symmetric and hence potential,
$$
A_{abc}= \frac{\p^3
\Phi^W}{\p t^a_W\p t^b_W\p t^c_W},
$$
for scalar function $\Phi^W(t_W)$ on $\cM$.
 The above quadratic
equations for $A_{ab}^c$ mean that the  latter
must be a solution of the WDVV equations. It is not hard to check
(cf.\ \cite{Ba2})
that the equation,
$$
\frac{\p\Psi^W(t_W, \hbar)}{\p \hbar}=-\hbar^{-1} \sum_c E^c(t_W)
\frac{\p\Psi^W(t_W, \hbar)}{\p t^c_W},
$$
holds for some vector field $E$ on $\cM$. Differentiating the
defining equation of the metric
with respect to $\hbar$, one
immediately concludes that $E$ is a conformal Killing vector with
respect to the metric $g$,
$$
{\cal L}_E g= (2-n)g.
$$
 Doing the same to the equation
defining the structure
functions, one proves the homogeneity,
$$
{\cal L}_E (\circ) =  \circ,
$$
of the product in $\cT_{\cM}$.

\sip

The basis $\{\Delta_a\}$ may be chosen in such a way that
$\Delta_0$ is the unit of the supercommutative algebra
$H^*(X,\C)$. We denote by $1$ the canonical lift of $\Delta_0$
into $\fg$.
The map of pointed dg manifolds, $\Gamma: (M,*,0) \rar (\fg, 0,
Q_{\fg})$, can always be normalized in such a way that \cite{Me2}
$$
\frac{\p \Gamma}{\p t^0_W} =1.
$$
This immediately implies,
$$
\frac{\p\Psi^W}{\p t^0_W}=  \hbar^{-1}\Psi^W,
$$
which in turn implies that the vector field $e=\p/\p t^0_W$ is a
unit with respect to the multiplication $\circ$.

\sip

Thus applying the remarkable Barannikov's construction \cite{Ba,Ba2}
to the deformation theory of K\"ahler forms, one gets
 a family, $(\Phi^W, g,E,e )$, of Frobenius manifold structures
on  $\cM\simeq H^*(X,\C)$ parameterized by isotropic filtrations
in $H^*(X,\wedge^*T_X)$ which are complementary to the standard
Hodge one.


\bip
\bip

\bip

\begin{center}
{\bf 5. Dual torus fibrations}
\end{center}

\bip

\sip

{\bf 5.1. Monge-Amp\`ere manifolds} \cite{H,KS}. Such a manifold
is a triple, $(Y,g,\nabla)$, consisting of a smooth manifold $Y$, a
smooth Riemannian metric $g$, and a flat torsion-free affine
connection $\nabla$ such that
\Bi
\item[(i)] the metric is potential in the sense that its
coefficients in a local affine coordinate system $\{x^i,
i=1,\ldots, \dim Y\}$ are given by
$$
g_{ij}= \frac{\p^2 K}{\p x^i \p x^j},
$$
for some smooth (convex) function $K$;
\item[(ii)] The  Monge-Amp\`ere equation,
$$
\det \left(\frac{\p^2 K}{\p x^i \p x^j}\right) = {\mathrm const},
$$
is satisfied.
\Ei

Due to the convexity of $K$ the system of equations,
$$
\frac{\p \hat{x}_i(x^j)}{\p x^j} = g_{ij},
$$
can solved in any local affine coordinate system $\{x^i\}$ on
$(Y,g,\nabla)$. The resulting functions $\hat{x}_i=\hat{x}_i(x^j)$ give rise
to a new local coordinate system on $Y$ and hence to a new
torsion-free affine connection, $\hat{\nabla}$, on $Y$ which has
$\{y_i\}$ as local affine coordinates.
Moreover, potentiality of the
metric and the Monge-Amp\`ere equation imply,
$$
g^{ij} =\frac{\p^2 \hat{K}}{\p y_i \p y_j},
$$
for some smooth function $\hat{K}$ satisfying
$$
\det \left(\frac{\p^2 \hat{K}}{\p y_i \p y_j}\right) = {\mathrm
const}^{-1}.
$$
Thus we have the following
\bip

{\bf 5.1.1. Proposition} \cite{H,KS,L}. {\em For any Monge-Amp\`ere
manifold $(Y,g,\nabla)$ there is canonically associated dual
Monge-Amp\`ere manifold $(\hat{Y}, \hat{g}, \hat{\nabla})$ such
that $(Y,g)=(\hat{Y},\hat{g})$ as Riemannian manifolds, while the
local systems $(T_Y, \nabla)$ and  $(T_{\hat{Y}}, \hat{\nabla})$
are dual to each other.
}

\bip

\sip

{\bf 5.2. Dual torus fibrations}
\cite{H,KS,L,SYZ}. Let $(Y,g,\nabla)$ be a Monge-Amp\`ere
manifold. The flat connection produces a natural splitting,
$T_{TY}= \pi^*(T_Y) \oplus \pi^*(T_Y)$, where $TY$ is the total space of
the tangent bundle to $Y$.  The almost complex
structure defined with respect to this splitting as $J:(v_1, v_2)
\rar (-v_2, v_1)$ is integrable and hence makes  $TY$ into a complex manifold.
Using the same splitting, one introduces a Riemannian metric on
$TY$,
$g_{TY}= \pi^*g \oplus \pi^*g$, which is K\"ahler with respect to
$J$. Its potential turns out to be $\pi^*(K)$ so that
the Monge-Amp\`ere equation for $K$ turns into the Ricci flatness of $g_{TY}$.
The net result is that the total space of the tangent bundle to a
Monge-Amp\`ere manifold is canonically a (non-compact) Calabi-Yau
manifold.

\sip

If the holonomy of the flat connection $\nabla$ is contained in
$SL(n,\Z)$, then one can choose a $\nabla$-parallel lattice
$T_Y^{\Z}$ and take the quotient, $X= TY/T_Y^{\Z}$. The result is
a Calabi-Yau manifold $X$ which is a torus fibration over $Y$
with typical fiber $T^n\simeq T_{Y,y}/T_{Y,y}^{\Z}$, $n=\dim Y$.
Applying the same sequence of constructions to the dual
Monge-Amp\`ere manifold $(\hat{Y}, \hat{g}, \hat{\nabla})$,
one gets another Calabi-Yau manifold $\hat{X}$ which is also a
torus fibration over $Y$. Now the  typical fiber is $\hat{T}^n \simeq
T_{Y,y}/ \hat{T}_{Y,y}^{\Z}$, where  $\hat{T}_{Y,y}^{\Z}$ is the
lattice in $T_{Y,y}$ which is dual to $T_{Y,y}^{\Z}$ with respect
to the metric $g$. If $Y$ is compact, it was argued by many authors
that the associated pair of Calabi-Yau manifolds, $X$ and
$\hat{X}$, is a {\em mirror pair} \cite{KS,L,SYZ}.

\bip

{\bf 5.3. Mirror symmetry between semi-infinite
variations of Hodge structure.} With any compact Calabi-Yau
manifold $X$ one can associate two models:
\Bi
\item semi-infinite $B$-variations of Hodge structure,
$$
 {\mathsf V}{\mathsf H}{\mathsf S}^B(X):
 (H^*(X, \wedge^* T_X),0) \lon Gr_{\frac{\infty}{2}}(H^*(X,\C)
[[\hbar,\hbar^{-1}]])(\f_0)
$$
 over the
extended moduli space of complex
structures, and a family, $\{\Phi^W_{\mathrm complex}(X)\}$, of
Frobenius manifold structures parameterized
by isotropic increasing filtrations, $W$, in $H^*(X,\C)$ which are
complementary to the Hodge one (this is the original Barannikov
construction, see \cite{Ba});

\item semi-infinite $A$-variations of Hodge structure,
$$
 {\mathsf V}{\mathsf H}{\mathsf S}^A(X): (H^*(X, \C),0) \lon
 Gr_{\frac{\infty}{2}}(H^*(X,\wedge^* T_X)
[[\hbar,\hbar^{-1}]])(\f_0)
$$
 over the
extended moduli space of K\"ahler forms,
and a family, $\{\Phi^{W'}_{\mathrm sympl}(X)\}$, of
Frobenius manifold structures parameterized
by isotropic increasing filtrations, $W'$, in $H^*(X,\wedge^* T_X)$ which are
complementary to the Hodge one (the symplectic version of the Barannikov
construction, see Sect.\ 4).
\Ei

{\bf 5.3.1. Theorem.} {\em
Let $X$ and $\hat{X}$ be dual
torus fibrations over a compact Monge-Amp\`ere manifold $Y$.
Then
$$
 {\mathsf V}{\mathsf H}{\mathsf S}^A(X) =  {\mathsf V}{\mathsf H}{\mathsf S}^B(\hat{X}), \ \ \
  {\mathsf V}{\mathsf H}{\mathsf S}^B(X) =  {\mathsf V}{\mathsf H}{\mathsf S}^A(\hat{X}).
$$
In particular,
$$
\Phi^W_{\mathrm complex}(X) = \Phi^{\hat{W}}_{\mathrm
sympl}(\hat{X}),
\ \ \ \Phi^{W'}_{\mathrm symp}(X) = \Phi^{\hat{W}'}_{\mathrm complex}(\hat{X})
$$
for appropriately related  filtrations $(W,\hat{W})$ and
$(W', \hat{W}')$.
}
\bip

\Proof  It is enough to study $T^n$-invariant subsheaves of
the sheaves of differential forms and polivector fields on both
$X$ and $\hat{X}$. Which are  related to each other via the
following composition of isomorphisms of $\f_Y$-modules,
$$
\phi: (\wedge^* T_X \ot \Omega^{0,*}_X)_{T^n}
\stackrel{i_i}{\lon} \wedge^*T_Y \ot \Omega^*_Y\ot \C
\stackrel{i_2}{\lon} \Omega^*_Y \ot \Omega^*_Y\ot \C
\stackrel{i_3}{\lon}
(\Omega^{*,*}_{\hat{X}})_{\hat{T}^n},
$$
and a similar one with $X$ and $\hat{X}$ exchanged. Here we implicitly
used the fact that $Y$ can be canonically embedded into both $X$ and $\hat{X}$ as
 the zero section. The middle
isomorphism,
$$
\wedge^*T_Y \ot \Omega^*_Y\ot \C
\stackrel{i_2}{\lon} \Omega^*_Y \ot \Omega^*_Y\ot \C
$$
comes from the ``lowering of indices'' map induced by the metric
$g$ on the first tensor factor and the identity maps on the other two factors.
In   corresponding local affine coordinate systems  this map is given on
generators by,
$$
\phi: \frac{\p}{\p z^i} \stackrel{i_1}{\lon}  \frac{\p}{\p x^i}
\stackrel{i_2}{\lon} \sum_j
g_{ij}dx^j = d\hat{x}_i
\stackrel{i_3}{\lon}
d\hat{z}_i,
$$
and
$$
\phi: d\bar{z}^i \stackrel{i_1}{\lon}
 dx^i \stackrel{i_2}{\lon} dx^i = \sum_j g^{ij}d\hat{x}_j \stackrel{i_3}{\lon}
\sum_j g^{ij} d\overline{\hat{z}}_j.
$$

Then we have, for example,
\Beqrn
\phi(\bp) &=& \phi (\sum_i d\bar{z}^i \frac{\p}{\p \bar{z}^i}) \\
&=&  \sum_{ij} g^{ij}  d\overline{\hat{z}}_j \frac{\p}{\p x^i}  \\
&=& \sum_j d\overline{\hat{z}}_j \frac{\p}{\p
{\hat{x}}_j}\\
&=& \sum_j d\overline{\hat{z}}_j \frac{\p}{\p
\overline{\hat{z}}_j}\\
&=& \bp,
\Eeqrn
and
\Beqrn
\phi(Q) &=& \phi\left(\sum_{ij} \psi_i \left(g^{ij} \frac{\p }{\p
\bar{z}^j} - \sum_k d\bar{z}^k \frac{\p g^{ij}}{\p \bar{z}^k}\frac{\p }{\p
d\bar{z}^j}\right) \right)\\
&=&\sum_{ij} dz_i \left(g^{ij} \frac{\p }{\p
x^j} + \sum_k dx^k \frac{\p g^{ij}}{\p x^k}\frac{\p }{\p
dx^j} \right)\\
&=&\sum_{ij} dz_i \left(g^{ij} \sum_k\left(\frac{\p x_k}{\p
x^j}\frac{\p}{\p x_k} +
\frac{\p dx_k}{\p x^j}\frac{\p}{\p dx_k}\right)  +
\sum_k dx^k \frac{\p g^{ij}}{\p x^k}\frac{\p }{\p
dx^j} \right)\\
&=& \sum_i dz_i \frac{\p}{\p x_i} \\
&=& \sum_i dz_i \frac{\p}{\p z_i}\\
&=& \p.
\Eeqrn
Analogously one checks that the map $\phi$ sends the Schouten
brackets into the Lie brackets defined in Sect.\ 2.6. As a
result,
the map $\phi$ canonically identifies the pair consisting of the dg Lie
algebra,
$$
\left(\bigoplus_{i=0}^{2n} \fg^i,\   \fg^i= \bigoplus_{p+q=i}
\Gamma(M, \wedge^p T_X\ot \Omega^{0,q}_{X})_{T^n},\ [\ \bullet\ ]_{\mathrm Sch},
\bp
\right),
$$
and its module,
$$
\left(\bigoplus_{i=0}^{2n} \fm^i,\   \fm^i= \bigoplus_{p+q=i}
\Gamma(X, \Omega^{p,q}_{X})_{T^n},\ \bullet\ ,\ \circ \ ,
d=\p+\bar{\p}\right),
$$
with the pair consisting of the dg Lie algebra
$$
\left(\bigoplus_{i=0}^{2n} \fg^i,\   \fg^i= \bigoplus_{p+q=i}
\Gamma(\hat{X}, \Omega^{p,q}_{\hat{X}})_{\hat{T}^n}\ ,\ [\ \bullet\ ]_{\omega},
\
\bar{\p}\right),
$$
and its module
$$
\left(\bigoplus_{i=0}^{2n} \fm^i,\ \fm^i= \bigoplus_{p+q=i}
\Gamma(\hat{X}, \Lambda^p T_{\hat{X}}\ot \Omega^{0,q}_{\hat{X}})_{\hat{T}^n}\ ,\ \bullet\ , \ \circ
\ , \bar{\p}+Q\right).
$$
The same statement, but with $X$ and $\hat{X}$ interchanged, is
also true. All the above claims follow immediately.
\hfill $\Box$

\bip

\pagebreak

\begin{center}
{\bf 6. Semi-infinite $\mathsf VHS$ in dGBV algebras}
\end{center}

\bip

\sip

{\bf 6.1. Differential Gerstenhaber-Batalin-Vilkovosky algebras.}
Such an algebra is a quadriple $(A, \circ, d, \Delta)$, where $(A,
\circ)$ is a unital supercommutative algebra over a field $k$, and
$(d,\Delta)$ is a pair of supercommuting  odd derivations of $(A,\circ)$ of order $1$ and $2$
respectively which satisfy $d^2=\Delta^2=0$.

\sip

Equivalently, a dGBV algebra is  a differential
supercommutative algebra with unit, $(A,\circ, d)$, plus an odd linear map
$\Delta: A \rar A$ satisfying
\Bi
\item[(i)] $\Delta^2=0$,
\item[(ii)] $d\Delta + \Delta d =0$,
\item[(iii)] and, for any $a,b,c\in A$,
\Beqrn
\Delta(a\circ b \circ c) & = & \Delta(a\circ b) \circ c +
(-1)^{\tl{b}(\tl{a}+1)} b\circ \Delta(a\circ c) +
(-1)^{\tl{a}}a\circ \Delta(b\circ c)\\
&& - \Delta(a) \circ b \circ c - (-1)^{\tl{a}} a\circ
\Delta(b)\circ c - (-1)^{(\tl{a}+\tl{b})} a\circ b \circ
\Delta(c).
\Eeqrn
 Note that $\Delta(1)=0$.
\Ei

\sip

It is not hard to check (see, e.g.\ \cite{Ma})  that the linear map
$$
\Ba{rccc}
[\ \bullet \ ] & : A\ot A & \lon & A \\
& a \ot b & \lon & [a\bullet b]:= (-1)^{\tl{a}} \Delta(a\circ b) -
(-1)^{\tl{a}} \Delta(a)\circ b - a\circ \Delta(b)
\Ea
$$
makes $A$ into a Lie superalgebra. Moreover, both the triples $(A, [\
\bullet\ ], d)$ and $(A, [\ \bullet \ ], \Delta)$ are differential
Lie superalgebras, and the following odd Poisson identity,
$$
[a\bullet (b\circ c)] = [a\bullet b]\circ c +
(-1)^{\tl{a}(\tl{b}+1)}b\circ [a\bullet b],
$$
holds for any $a,b\in A$.

\sip

\bip

{\bf 6.2. A dGBV algebra as a $\fg$-module.}
For any $a\in A$, we denote
$$
\Ba{rccc}
l_a: & A& \lon &   A \\
& b & \lon & l_a(b):= a\circ b.
\Ea
$$
The parity of this linear map is equal to the parity of $a$.

\bip

{\bf 6.2.1. Proposition.} {\em Let $(A,\circ, d , \Delta)$ be a dGBV
algebra. Then the linear map,
$$
\Ba{rccc}
\bullet & : A\ot A & \lon A \\
& a \ot b & \lon  & a\bullet b:= -[l_a, \Delta]b,
\Ea
$$
makes the triple $\fm =( A, \bullet, d+\Delta)$ into a
differential module over the differential Lie algebra $\fg=(A,
[\ \bullet\ ], d)$.
}

\bip

\Proof
Let us first show that $(A,\bullet)$ is a module over the Lie
algebra $(A,[\ \bullet\ ])$, that is,
$$
a_1\bullet a_2 \bullet b -
(-1)^{(\tl{a}_1+1)(\tl{a}_2+1)}a_2\bullet a_1\bullet b =
[a_1\bullet a_2] \bullet b,
$$
holds for all $a_1,a_2,b\in A$.

\sip

Using the definition and the nilpotency of $\Delta$, we have
\Beqrn
{\mathrm l.h.s.} &=& a_1\circ \Delta(a_2 \circ \Delta b) -
2(-1)^{\tl{a}_1} \Delta(a_1\circ a_2 \circ \Delta b)
 + (-1)^{\tl{a}_1 +
\tl{a}_2} \Delta(a_1\circ \Delta(a_2 \circ b))\\
&&+ (-1)^{\tl{a}_1\tla_2 + \tla_1 + \tla_2} a_2\circ \Delta(a_1 \circ \Delta b)
+ (-1)^{\tla_1\tla_2} \Delta(a_2\circ \Delta(a_1 \circ b)).
\Eeqrn
On the other hand, using  definitions only, we have
\Beqrn
{\mathrm r.h.s.} &=& - [l_{[a_1\bullet a_2]}, \Delta] b \\
&=& - (-1)^{\tla_1} \Delta(a_1\circ a_2) \circ \Delta b +
(-1)^{\tl{a}_1} \Delta(a_1)\circ a_2\circ \Delta b + a_1\circ
\Delta(a_2) \circ \Delta b \\
&&- (-1)^{\tla_2}\Delta\left( \Delta(a_1\circ a_2)\circ b -
\Delta(a_1)\circ a_2\circ b - (-1)^{\tla_1} a_1\circ
\Delta(a_2)\circ b\right).
\Eeqrn
As $\Delta \Delta(a_1\circ a_2\circ b)$ vanishes identically,

\bip

\noindent\mbox{$\Delta\left( \Delta(a_1\circ a_2)\circ b -
\Delta(a_1)\circ a_2\circ b - (-1)^{\tla_1} a_1\circ
\Delta(a_2)\circ b\right)=$} \newline
\mbox{\ \ \ \ \ \ \ \ \ \ \ \
$- \Delta\left( (-1)^{\tla_2(\tla_1+1)}a_2\circ
\Delta(a_1\circ b) + (-1)^{\tla_1}a_1\circ \Delta(a_2\circ b) -
(-1)^{\tla_1+\tla_2} a_1\circ a_2\circ \Delta b\right)$},

\bip

\noindent so that using the decomposition formula for
$\Delta(a_1\circ a_2\circ \Delta b)$ one gets eventually  the
desired equality
$$
{\mathrm r.h.s.} = {\mathrm l.h.s.}
$$

Next,
\Beqrn
(d+\Delta) (a_1\bullet b) &=& (d+\Delta)(-a_1\circ \Delta b +
(-1)^{\tla_1} \Delta(a_1\circ b))\\
&=& -da_1 \circ \Delta b + (-1)^{\tla_1}a_1 \circ \Delta(d+\Delta)b
- (-1)^{\tla_1}\Delta(da_1\circ b)\\
&& -\Delta(a_1 \circ db) - \Delta(a_1 \circ \Delta b)\\
&=& da_1\bullet b - (-1)^{\tla_1} a_1 \bullet (d+\Delta)b,
\Eeqrn
confirming the consistency  with the differentials.
\hfill $\Box$

\sip

\bip

{\bf 6.3. Local systems from dGBV-algebras.} In the notions of
 Proposition~6.2.1, we reinterpret the algebra structure in $A$ as  an even
 map,
$$
\Ba{rccc}
\circ: & \fg \ot \fm & \lon & \fm\\
& a\ot b & \lon & a \circ b.
\Ea
$$

\bip

{\bf 6.3.1. Proposition.} {\em For any $a_1,a_2,b\in A$, we have
$$
(d+\Delta)(a_1 \circ b)= (da_1)\circ b + (-1)^{\tla_1}a_1\circ
(d+\Delta)b + (-1)^{\tla_1}a_1\bullet b,
$$
and
$$
a_1\circ a_2 \bullet b - (-1)^{\tla_1(\tla_2 + 1)}a_2\bullet a_1
\circ b =  -(-1)^{\tla_2} [a_1\bullet a_2] \circ b.
$$
}

\bip

\Proof The first equality is obvious. We check the second one,
\Beqrn
{\mathrm l.h.s.} &=& -a_1 \circ a_2 \circ \Delta b +
(-1)^{\tla_2}a_1\circ\Delta(a_2\circ b) + (-1)^{\tla_1(\tla_2 +
1)}a_2\circ (a_1\circ b)\\
&& - (-1)^{\tla_1+\tla_2}\Delta(a_1\circ a_2 \circ b)\\
&=& -(-1)^{\tla_1+\tla_2}\Delta(a_1\circ a_2) \circ b +
(-1)^{\tla_1+\tla_2}\Delta(a_1)\circ a_2 \circ b +(-1)^{\tla_2}a_1\circ \Delta(a_2) \circ
b\\
&&= {\mathrm r.h.s.}.
\Eeqrn
\hfill $\Box$

\bip

By Proposition~3.2.1, the dGBV algebra $(A, \circ, d, \Delta)$
gives rise to a flat vector bundle, $(E_{\fm}, \nabla)$, over the
derived moduli space $\cM$ representing the deformation functor $\Def_{\fg}$.
The typical fibre of $\pi: E_{\fm}\rar \cM$ is
isomorphic to the cohomology group $H(A, d+\Delta)$.

\sip

\bip

{\bf 6.4. Semi-infinite variations of Hodge structure in $A$.}
There is an obvious formal parameter extension of all the results in
Sections 6.2 and 6.3. More precisely,
for any dGBV algebra $(A,\circ, d,\Delta)$, the pair,
$(\fg, \fm_{\hbar})$,  consisting of the  differential Lie superalgebra
$$
\fg = (A,\ [\ \bullet\ ], d),
$$
and its differential module,
$$
\fm_{\hbar}:= (A[[\hbar, \hbar^{-1}]], \bullet_{\hbar}:=\bullet,
\circ_{\hbar}:=\hbar^{-1}\circ, d+\hbar \Delta),
$$
do satisfy the conditions of
Proposition~3.2.1 and hence give rise to a local system,
$(E_{{\hbar}}, \nabla)$ over the derived moduli space $\cM$.

\sip

\bip

{\bf Proposition 6.4.1.} {\em Assume that the dGBV algebra $A$
satisfies,
$$
\Img d \cap \Ker \Delta = \Img \Delta \cap \Ker d = \Img d \cap
\Img \Delta,
$$
then
\Bi
\item[(i)] the deformation functor $\Df_{\fg}$ is non-obstructed
so that
$
(\cM, *) \simeq (H(A,d), 0)
$
as formal pointed supermanifolds;
\item[(ii)] the typical fibre of the associated bundle, $\pi:
E_{{\hbar}}\rar \cM$, is isomorphic to
$H(A,d)[[\hbar,\hbar^{-1}]]$.
\Ei
}
\bip

\Proof It follows from \cite{Ma} that $\eth=0$ and that
$H(A,d)=H(A,\Delta)= \Ker d \cap \Delta/ \Img d\Delta$.
\hfill $\Box$

\bip

Assuming further  that the dGBV algebra is $\Z\times \Z$-graded (or $\Z$-graded),
$A=\oplus_{i,j}A^{i,j}$, with $\circ$ of bidegree $(0,0)$,  $d$ of bidegree
$(0,1)$, $\Delta$ of bidegree $(-1,0)$, and $\dim H(A,d)<\infty$,
 we can use Barannikov's technique \cite{Ba} to define
 semi-infinite variations of the Hodge structure in $H(A,d)$, and
 then
  construct a family of pencils of flat connections on $\cM$, $\nabla^W$,
 depending on the choice of increasing filtrations, $W$, which are
 complementary to the Hodge one.

 \sip

 Moreover, if $(A,\circ, d, \nabla)$ admits a linear map
(called  an {\em integral}\, in \cite{Ma}),
$$
\int: A \lon k
$$
satisfying
$$
\int (da) \circ b = -(-1)^{\tl{a}+1}\int a\circ db, \ \ \ \ \forall
a,b\in A,
$$
$$
\int (\Delta a)\circ b = (-1)^{\tl{a}}\int a\circ \Delta b, \ \ \ \ \forall
a,b\in A \\
$$
and inducing a non-degenerate metric, $g([a], [b])= \int a\circ
b$, on $H(A,d)$, then one can make one step further and construct
a family of solutions, $\{\Phi^{W}\}$, of WDVV equations on $\cM$
which are parameterized by the $g$-isotropic subclass of the above class of filtrations.
This links Barannikov's construction with the one
studied in \cite{Ma}.

\bip

{\bf 6.5. Examples.} (i) Let $(M, \omega)$ be a Lefschetz
symplectic manifold, then the data
$$
\left(A=\Gamma(M, \Omega^*_M), \circ={\mathrm wedge\ product}, d,
\Delta=[d,\Lambda_{\omega}]\right)
$$
where $\Lambda_{\omega}: \Omega^*_M \rar \Omega^{*-2}$ is the map
of contraction with $\omega^{-1}$ is a dGBV algebra satisfying the
conditions of Proposition~6.4.1 and admitting an integral.

\sip

(ii) If $(M,\omega)$ is a K\"ahler manifold, then the following
variant of the above data,
$$
\left(A=\Gamma(M, \Omega^{*,*}_M), \circ={\mathrm wedge\ product}, \bp,
\Delta=[\bp,\Lambda_{\omega}]\right)
$$
is a dGBV algebra satisfying all the conditions of subsection~6.4 necessary for
constructing  solutions, $\{\Phi^W\}$, of WDVV equations from
semi-infinite variations of the Hodge structure. If $M$ is
Calabi-Yau, this construction is equivalent to the one discussed in
Section 4, but in general it is different.

\sip

(iii) If $M$ is a Calabi-Yau manifold with a global holomorphic
volume form $\Omega$, then
$$
\left(A=\Gamma(M, \wedge^* T_M \ot \Omega^{0,*}_M), \circ={\mathrm wedge\ products}, \bp,
\Delta=i_{\Omega}^{-1} \p i_{\Omega}\right)
$$
with $i_{\Omega}$ being the isomorphism, $\wedge^pT_M\ot \Omega^{0,q}_M
\rar \Omega^{n-p,q}_M$, induced by the volume form, is a dGBV
algebra satisfying all the conditions of subsection 6.4. It is easy to check that the
resulting construction of Frobenius
manifolds is equivalent to the original Barannikov's one
\cite{Ba}.

\vse

{\small

{\small
\begin{tabular}{l}
Department of Mathematics\\
University of Glasgow\\
15 University Gardens \\
 Glasgow G12 8QW, UK
\end{tabular}
}

\end{document}